\documentclass[12pt,reqno]{amsart}
\pdfoutput=1
\usepackage[headings]{fullpage}
\usepackage{amsfonts}
\usepackage{amsthm}
\usepackage{amssymb}
\usepackage{epstopdf}
\usepackage{graphicx}
\usepackage{texdraw}
\usepackage{url}
\usepackage{xfrac}
\usepackage{enumitem}
\usepackage{float}   
\usepackage{mathtools}  
\usepackage[bookmarks=true,%
    colorlinks=true,%
    linkcolor=blue,%
    citecolor=blue,%
    filecolor=blue,%
    menucolor=blue,%
    urlcolor=blue,%
    breaklinks=true]{hyperref}

\newtheorem{theorem}{Theorem}[section]

\theoremstyle{definition}

\def\BZ{\mathbb Z}

\def\BQ{\mathbb Q}

\def\BC{\mathbb C}

\def\BH{\mathbb H}

\def\calT{\mathcal T}

\def\Ga{\Gamma}

\def\longto{\longrightarrow}

\def\PSL{\mathrm{PSL}}

\def\pt{\partial}
\def\pr{\mathrm{par}}
\def\be{\begin{equation}}
\def\ee{\end{equation}}

\begin{document}


\title[Factorization of polynomials in hyperbolic geometry and dynamics]{
  Factorization of polynomials in hyperbolic geometry and dynamics}
\author{Michael Filaseta}
\address{Department of Mathematics \\
         University of South Carolina \\
         Columbia, SC 29208, USA \newline 
         {\tt \url{http://www.math.sc.edu/~filaseta }}}
\email{filaseta@mailbox.sc.edu}
\author{Stavros Garoufalidis}
\address{
  International Center for Mathematics, Department of Mathematics \\
  Southern University of Science and Technology \\
  Shenzhen, China \newline
  {\tt \url{http://people.mpim-bonn.mpg.de/stavros}}}
\email{stavros@mpim-bonn.mpg.de}
\date{17 September 2022}


\begin{abstract}
Using factorization theorems for sparse polynomials, 
we compute the trace field of Dehn fillings of the Whitehead link, and
(assuming Lehmer's Conjecture) the minimal polynomial of the small 
dilatation pseudo-Anosov maps and the trace field of fillings of the
figure-8 knot.  These results depend on the degrees of the trace fields over
$\mathbb Q$ being sufficiently large.
\end{abstract}

\maketitle

{\footnotesize
\tableofcontents
}


\section{Introduction}
\label{sec.intro}

\subsection{Factorization problems in dynamics and 
hyperbolic geometry}
\label{sub.problems}

Two problems that one encounters in dynamics and hyperbolic geometry
are about number fields and their defining polynomials, namely the
computation 
\begin{itemize}
\item
of the minimal polynomial of the largest positive eigenvalue of
a Perron-Frobenius matrix, and of 
\item
of the trace fields of hyperbolic Dehn-fillings.
\end{itemize}
In both problems the number fields are generated by a solution of a one-variable
specialization $p(x^m,x^n)$ of a polynomial $p(x,y) \in \BZ[x,y]$ for some integers
$m$ and $n$. This leads to the question of the factorization of $p(x^m,x^n)$
for $p(x,y) \in \BZ[x,y]$. It is customary to factorize a polynomial with
integer coefficients into its cyclotomic part, its reciprocal non-cyclotomic
part and its non-reciprocal part, defined in detail in Section~\ref{sec.factoring1}.
Thus, the above problem reduces to the factorization of the cyclotomic,
reciprocal and non-reciprocal parts of $p(x^m,x^n)$. The cyclotomic part is
well-understood; see for example the work of Granville and
Rudnick~\cite[Section 3]{Granville}. The non-reciprocal part is also well-understood
by the work of Schinzel~\cite{schinzel-I}; see also~\cite{filaseta97}
and~\cite[Theorem 1.1]{FMV}. 
On the other hand, the factorization of the reciprocal part of a polynomial is
not generally known, abstractly or concretely.

In special circumstances regarding the trace fields of hyperbolic Dehn-fillings,
one can bypass the detailed problem of the factorization of $p(x^m,x^n)$,
and assuming Lehmer's Conjecture, obtain a qualitative result asserting that
each irreducible factor has degree bounded below by $C \max\{|m|,|n|\}$ for
some constant that depends on $p$; see~\cite{Ga-Jeon}. 
Furthermore, such results are also possible unconditionally by using the work
of Dimitrov~\cite{Dimitrov}, though the special circumstances become more restrictive.

In the current paper, we will use the above described factorization methods
in two sample problems.  In a subsequent paper \cite{Filaseta22}, an approach is developed 
to obtain information about the factorization of polynomials of the form $f_{0}(x) + f_{1}(x) x^{n} + \cdots + f_{r-1}(x) x^{(r-1)n} + f_{r}(x) x^{rn}$,
replacing the role of Theorem~\ref{ydeg1} below and leading to weaker results
than those in this paper but for polynomials associated with more general trace fields 
of hyperbolic Dehn-fillings. 

\subsection{The trace field of hyperbolic Dehn fillings}

The fundamental group of an oriented hyperbolic 3-manifold $M$ embeds as a 
discrete subgroup $\Ga$ of $\PSL(2,\BC)$ and the embedding is well-defined up 
to conjugation \cite{Th}. Discreteness implies that the subfield of the
complex numbers generated by the traces of the elements of $\Ga$ is a 
number field, the so-called {\em trace field} of $M$. For a detailed
discussion, see \cite[Thm.3.1.2]{MReid}.  

If the hyperbolic manifold $M$ has a cusp, i.e., a torus boundary component,
then all but finitely many fillings of a hyperbolic manifold are 
hyperbolic. This is the content of Thurston's hyperbolic Dehn filling theorem
\cite{Th,NZ}. Dehn fillings are parametrized by a rational number $m/n$.
A natural question is to ask how the trace field of the filled
manifold $M_{m/n}$ depends on the filling $m/n$.

Using ideal triangulations of hyperbolic manifolds and their
gluing equations and character varieties, one can reduce the problem of 
computing the trace field of an $M_{m/n}$ (for fixed $m$ and large $n$) to
the factorization of a {\em sparse polynomial} $p(x^{m},x^n)$ where $p(x,y) \in
\BZ[x,y]$ is a polynomial of positive $y$-degree. This is explained in 
detail in Section \ref{sec.dehn} of the appendix.

\begin{figure}[htpb]
\begin{center}
\includegraphics[height=0.07\textheight]{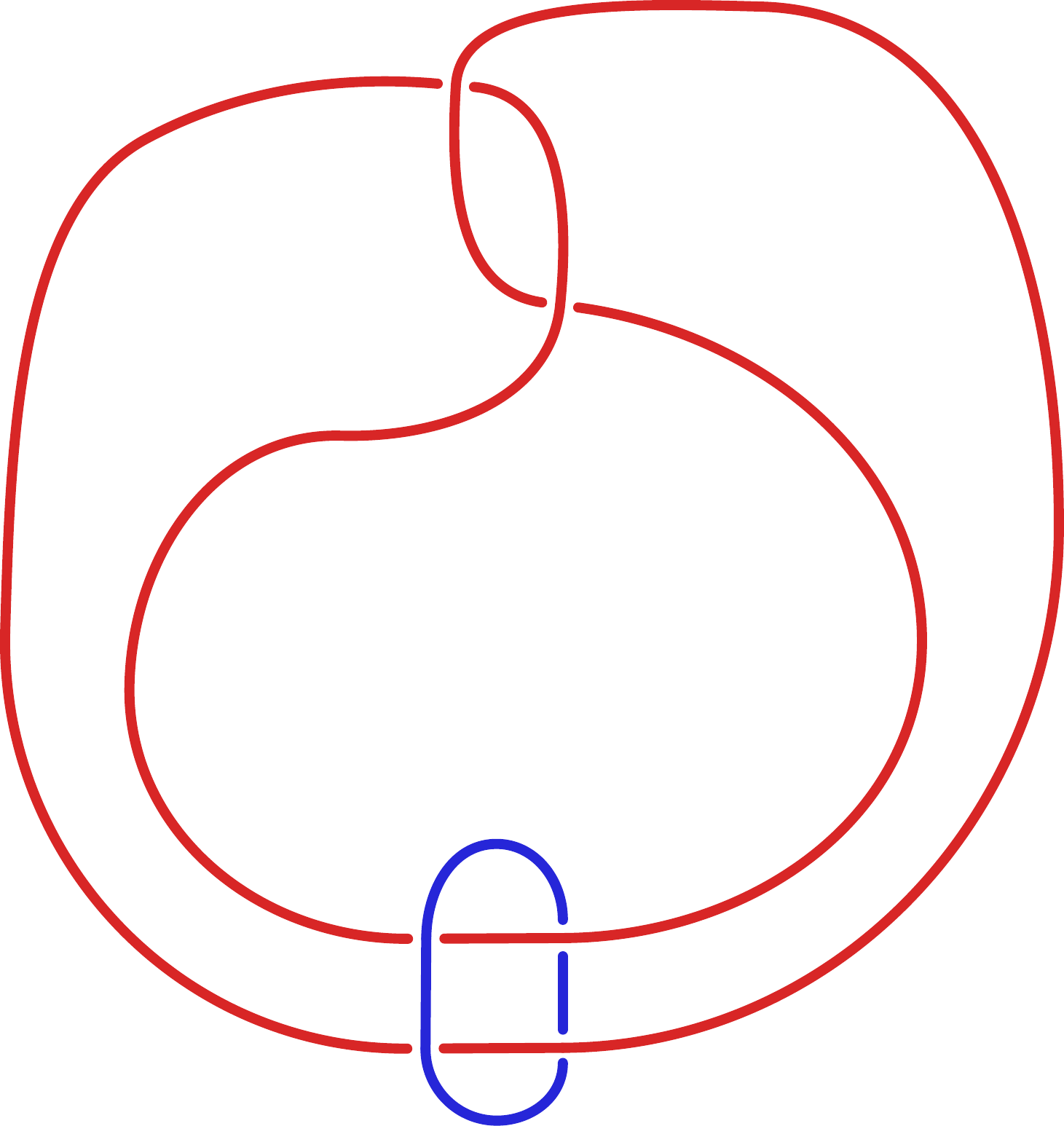}    
\caption{The Whitehead link.}
\label{f.2links}
\end{center}
\end{figure}

For example, consider the hyperbolic manifold $W$ with two cusps which is the 
complement of the Whitehead link shown in Figure \ref{f.2links}.

Neumann-Reid show 
in~\cite[Thm.~6.2, Eqn.~(6.7)]{NR} that the trace field of $W_{m/n}$ is given
by $\BQ(x-1/x)$ where $x$ is a suitable non-cyclotomic root of the polynomial
\begin{equation}
\label{eq.Fwhitehead}
F^W_{m,n}(x)=(x(x+1))^m x^{4n}-(x-1)^m.
\end{equation}
It is easy to see that $x^2+1$ is a root of $F^W_{m,n}(x)$.
Theorem~\ref{thm.whitehead} below proves that the quotient is an irreducible
polynomial in $\BZ[x]$ for every fixed $m$ and sufficiently large $n$.

A 2-variable polynomial of degree 1 in $y$ whose specialization
is~\eqref{eq.Fwhitehead} is given by
\[
P^W_m(x,y)=(x+1)^m y -(x-1)^m. 
\]
Section~\ref{sec.factoring1} deals with the non-reciprocal factors of
specializations of the above polynomial given by $y = x^{4n+m}$ 
where $n$ and $m$ are positive relatively prime integers with $m$ odd and $n$ large.
 
Our first result concerns fillings of the Whitehead link.

\begin{theorem}
\label{thm.whitehead}
Let $m$ be a fixed positive odd integer.  There is an integer $N = N(m)$ such
that if $n \in \mathbb Z$ with $n \ge N$ and $\gcd(n,m) = 1$, then 
$F^W_{m,n}(x)$ 
is $x^{2}+1$ times an irreducible polynomial. 
\end{theorem}

It is conceivable that $N(m) = 1$ for all $m$. Our approach provides an explicit
although impractical value for $N(m)$ for a given $m$, namely
\[
N(m) = \dfrac{5^{8 \binom{2m}{m} +8m - 7}}{2}- \dfrac{m}{4}.
\]
A corollary of the above theorem is the determination of the invariant trace field
$F_{W(m/n)}$ of the hyperbolic manifold $W(m/n)$ obtained by $m/n$ surgery
on the Whitehead link $W$. 

\begin{theorem}
\label{thm.whitehead2}
For $m$ a positive odd integer and $n$ an integer $> N(m)$, as defined above,
with $\gcd(n,m) = 1$, we have $F_{W(m/n)}=\BQ(z)$ where $z = x-1/x$ with $x$ a root
of 
\[
((x(x+1))^m x^{4n}-(x-1)^m)/(x^2+1).
\]
The degree of $F_{W(m/n)}$ is $2n+m-1$,
and a defining polynomial $T_{m,n}(z)$ of $F_{W(m/n)}$ over $\mathbb Q$ satisfies
linear recursions with respect to $m$ and $n$ with coefficients in $\BQ[z]$.
Specifically, we can take
\begin{equation}
  \label{Tmnrec}
\begin{aligned}
T _{m,n}(z) &= (z+2)\, T _{m-1,n}(z) - z \,T _{m-2,n}(z)  \quad \text{for $m \ge 2$}, \\
T _{m,n}(z) &= (z^2+2)\, T _{m,n-1}(z) - T _{m,n-2}(z)
\quad \text{for $n \ge 2$},
\end{aligned}
\end{equation}
with initial conditions
\be
  T_{0,0}(z) = 0, \quad T_{1,0}(z) = 1, \quad T_{0,1}(z) = z,
  \quad T_{1,1} = z^{2} + z + 1 .
\ee
\end{theorem}

\subsection{The minimal polynomial of small dilatation pseudo Anosov 
maps}
\label{sub.teichmuller}

Consider the 2-variable Laurent polynomial
\begin{equation}
\label{eq.teichmuller}
P^T(x,y) = y+y^{-1} -(x+x^{-1}+1).
\end{equation}
This is the Teichmuller polynomial of the simplest pseudo-Anosov braid
$\sigma_1 \sigma_2^{-1}$ studied in~\cite[p.552]{mcmullen2}. 
The specialization $P^T(x,x^g)$ appears in \cite[Thm.1.1]{mcmullen:entropy}
in connection with the minimum value of the spectral radius of all
reciprocal Perron-Frobenius $g \times g$ matrices. The specialization
$P^T_{a,b}(x):=P^T(x^a, x^b)$ for
positive integers $a$ and $b$ (where $a$ is within $2$ from a large integer $g$ and
$b$ is small) appears in relation to the smallest dilatation number
of pseudo-Anosov genus $g$ fiberings~\cite[Sec.1.1]{dunfield},
and also in~\cite{LT,hironaka}.

Clearing the denominators in $P^T(x,y)$ to obtain a polynomial in $\mathbb Z[x,y]$
gives
\[
x y P^T(x,y) = x y^{2} + x - x^{2} y - y - x y.
\]
We will be interested in the polynomial
\begin{equation}\label{GTdef}
G^{T}_{a,b}(x)
= \dfrac{x^{a} x^{b} P^T(x^{a},x^{b})}{x^{\min\{ a,b \}}}
= \dfrac{x^{a+2b} + x^{a} - x^{2a+b} - x^{b} - x^{a+b}}{x^{\min\{ a,b \}}}
\end{equation}
in $\mathbb Z[x]$ having constant term $\pm 1$.

Computations suggest that $G^T_{a,b}(x)$ has the following 
cyclotomic factors.
\begin{enumerate}[leftmargin=28pt]
\item[(I)]
  If $(a,b)$ is of the form $(a_{1}d_{1},b_{1}d_{1})$ where $d_{1} \in \mathbb Z^{+}$
  and $(a_{1},b_{1})$ is congruent to one of the pairs $(1, 0)$, $(3, \pm 2)$, and
  $(5,0)$ modulo $6$, then $G^T_{a,b}(x)$ is divisible by $\Phi_{6 d_{1}}(x)$.
\vskip 5pt \noindent
\item[(II)]
  If $(a,b)$ is of the form $(a_{2}d_{2},b_{2}d_{2})$ where $d_{2} \in \mathbb Z^{+}$
  and $(a_{2},b_{2})$ is congruent to one of the pairs
$(2, \pm 1)$, $(4,\pm 3)$, $(6,\pm 3)$ and $(8,\pm 1)$ modulo $10$, then
$G^T_{a,b}(x)$ is divisible by $\Phi_{10 d_{2}}(x)$.
\vskip 5pt \noindent
\item[(III)]
  If $(a,b)$ is of the form $(a_{3}d_{3},b_{3}d_{3})$ where $d_{3} \in \mathbb Z^{+}$
  and $(a_{3},b_{3})$ is congruent to one of the pairs $(3, \pm 2)$, $(4, \pm 3)$,
  $(8, \pm 3)$, and $(9,\pm 2)$ modulo $12$, then $G^T_{a,b}(x)$ is divisible by
  $\Phi_{12d_{3}}(x)$.
\end{enumerate}
(Here, $\Phi_n(x)$ is the $n$-th cyclotomic polynomial.) 
A given pair $(a,b)$ can occur in several cases. 
For example, if $(a,b) = (448,441)$, then (I) does not arise since $b$ is odd so
that $b_{1}$ cannot be $0$ or $\pm 2$ modulo $6$. However, each of (II) and (III)
can arise in two ways.  
For (II), we see that $(448,441) \equiv (8,1) \pmod{10}$ so taking $d_{2} = 1$ 
we obtain $\Phi_{10}(x)$ is a factor of $G^T_{a,b}(x)$; and
$(448,441) = (64 \cdot 7, 63 \cdot 7)$ where $(64,63) \equiv (4,3) \pmod{10}$ 
so taking $d_{2} = 7$ we obtain $\Phi_{70}(x)$ is a factor of $G^T_{a,b}(x)$. 
For (III), we see that $(448,441) \equiv (4,-3) \pmod{12}$ so taking $d_{3} = 1$ 
we obtain $\Phi_{12}(x)$ is a factor of $G^T_{a,b}(x)$; and
$(448,441) = (64 \cdot 7, 63 \cdot 7)$ where $(64,63) \equiv (4,3) \pmod{12}$ 
so taking $d_{3} = 7$ we obtain $\Phi_{84}(x)$ is a factor of $G^T_{a,b}(x)$. 
Thus, $G^T_{448,441}(x)$ has four cyclotomic factors. Indeed, a computation shows that
\be
G^T_{448,441}(x) = \Phi_{10}(x) \Phi_{70}(x) \cdot \Phi_{12}(x) \Phi_{84}(x)
\cdot(\text{non-cyclotomic irreducible polynomial}) .
\ee

Setting $C^T_{a,b}(x)$ to be the product of the distinct cyclotomic polynomials
arising from (I), (II) and (III), we will show that, for large $a/b$ or large
$b/a$, we have $G^T_{a,b}(x)$ is $C^T_{a,b}(x)$ times a non-cyclotomic
irreducible polynomial under the assumption of Lehmer's Conjecture on the
smallest Mahler measure of a non-cyclotomic irreducible polynomial with
integer coefficients. More precisely, we establish the following.

\begin{theorem}
\label{thm.teichmuller}
Let $a$ and $b$ be positive integers.  Then the following hold.
\begin{enumerate}[leftmargin=20pt]
\item[(a)]
The cyclotomic part of $G^T_{a,b}(x)$ is $C^T_{a,b}(x)$.
\item[(b)]
If $a$ and $b$ are distinct, then the polynomial
$G^T_{a,b}(x)$ is the product of irreducible reciprocal factors.
\item[(c)]
Assuming Lehmer's Conjecture, there exists an absolute constant $C$ such that if 
\[
\max\{ a/b, b/a \} > C,
\]
then the polynomial $G^T_{a,b}(x)/C^T_{a,b}(x)$ is irreducible. 
\end{enumerate}
\end{theorem}

In particular, if $g$ is a sufficiently large integer, then Lehmer's Conjecture
implies that $x^{g} P^T(x,x^g)/C^T_{1,g}(x)$ is an irreducible polynomial of
degree $2g-2\epsilon_g$ where $\epsilon_g=1$ if $6$ divides $g$ and $0$ otherwise.


\section{Factoring $f(x)x^n+g(x)$}
\label{sec.factoring1}

\subsection{Background on factorization}
\label{sub.background}

Let us begin by recalling some standard terminology regarding factorization. 
We define a non-zero polynomial $w(x) \in \mathbb R[x]$ as \textit{reciprocal}
if $w(x) = \pm x^{\deg w} w(1/x)$.  For example, $x+1$, $x - 1$ and 
$x^{3} - 3x^{2} + 3x -1$ are all reciprocal.  Observe that $w(0) \ne 0$ for a
reciprocal polynomial $w(x)$; to see this simply observe that if $x$ divides
$w(x)$, then the degree of $x^{\deg w} w(1/x)$ is less than the degree of
$w(x)$ so $w(x) \ne \pm x^{\deg w} w(1/x)$.  The terminology reciprocal is 
motivated by the fact that an (easily seen) equivalent definition for 
$w(x)$ being reciprocal is that $w(0) \ne 0$ and if $\alpha$ is a root of
$w(x)$ with multiplicity some positive integer $k$, then $1/\alpha$ is a root
of $w(x)$ with multiplicity $k$.  Observe that if $w(x)$ is a reciprocal
factor of a non-zero polynomial $F(x) \in \mathbb R[x]$, then $w(x)$ 
divides $x^{\deg F} F(1/x)$.  

Important examples of reciprocal polynomials are given by the classical cyclotomic 
polynomials, that is the irreducible factors of $x^{k}-1$ where
$k$ is a positive integer.  There is one irreducible factor $\Phi_k(x)$
of $x^{k}-1$ that does not appear as a factor of $x^{\ell}-1$ for $\ell$
a positive integer $< k$, and that irreducible factor is called the
$k^{\rm th}$ cyclotomic polynomial. It has $\zeta = e^{2 \pi i/k}$ as a root.  

Given a polynomial $F(x) \in \mathbb Z[x]$, we refer to the non-reciprocal 
part of $F(x)$ as the factor that remains after we remove each reciprocal factor
of $F(x)$ that is irreducible over the integers (so the gcd of the coefficients
is $1$) and that has a positive leading coefficient. Note that a non-zero
constant polynomial in $\mathbb Z[x]$ is reciprocal, and, as irreducibility
is over the integers, it is irreducible if and only if it is prime.
For example, the non-reciprocal part of 
\[
F(x) = (-2x+2)(x^6+x^5+x^4+3x^3+x^2+x+1)
\]
is obtained by first rewriting and completing the factorization to
\[
F(x) = -2(x-1)(x^3+x^{2}+1)(x^{3}+x+1)
\]
and then removing the irreducible reciprocal factors with a positive 
leading coefficient, namely $2$ and $x-1$. 
Thus, the non-reciprocal part of $F(x)$ in this example is
\[
-x^6-x^5-x^4-3x^3-x^2-x-1.
\]
As this example shows, the non-reciprocal part of $F(x)$ can itself be a 
reciprocal polynomial (but its
irreducible factors must be non-reciprocal).

For $w(x) = \sum_{j=0}^{r} a_{j} x^{j} \in \mathbb R[x]$, we use the notation
\[
\Vert w \Vert = \sqrt{a_{0}^{2} + a_{1}^{2} + \cdots + a_{r}^{2}}.
\]
Thus, $\Vert w \Vert$ may be viewed as a way for us to give some notion of 
a size to a polynomial.  

For $u(x)$ and $v(x)$ in $\mathbb Z[x]$, 
we define $\gcd{}_{\mathbb Z}(u(x),v(x))$
\index{$\gcd{}_{\mathbb Z}(u(x),v(x))$} as follows.
Write 
\[
u(x) = \pm \prod_{j=1}^{r} g_{j}(x)^{k_{j}} 
\quad \text{ and } \quad v(x) = \pm \prod_{j=1}^{r} g_{j}(x)^{\ell_{j}},
\]
where the $g_{j}(x)$ are distinct irreducible polynomials in $\mathbb Z[x]$ 
with positive leading coefficients and the exponents $k_{j}$ and $\ell_{j}$
are all non-negative integers. Then
\[
\gcd{}_{\mathbb Z}(u(x),v(x)) = \prod_{j=1}^{r} g_{j}(x)^{\min\{ k_{j}, \ell_{j} \}}.
\]
The example
\[
\gcd{}_{\mathbb Z}\big(-5(2x+2)^{2}(x^{2}-3)(x^{3}+x+1),-10(2x+2)^{3}(x^{2}-3)\big)
= 5(2x+2)^{2}(x^{2}-3)
\]
illustrates this definition.

Schinzel~\cite{schone,schinzel2}, capitalizing on an idea of
Ljunngren~\cite{wljun}, showed how one can obtain information on
the factorization of the non-reciprocal part of a polynomial.  
For the proof of Theorem \ref{thm.whitehead}, we make use of the following 
result from~\cite[p.635]{FFK}, which is a variation of a result that appears
in~\cite{schone}.

\begin{theorem}
\label{ydeg1}
Let $f(x)$ and $g(x)$ be in $\BZ[x]$ with $f(0) \ne 0$, $g(0) \ne 0$,
and $\gcd_{\BZ}(f(x),g(x)) = 1$. Let $r_{1}$ and $r_{2}$ denote the number
of non-zero terms in $f(x)$ and $g(x)$, respectively. If
$$
n \ge \max \bigg\{ 2 \times 5^{2N-1}, 
2 \max \big\{ \deg f, \deg g \big\} \bigg( 5^{N-1} + \dfrac{1}{4} \bigg) 
\bigg\}
$$ 
where
$$
N = 2\,\Vert f\Vert^{2} + 2\,\Vert g\Vert^{2} + 2r_{1} + 2r_{2} - 7,
$$
then the non-reciprocal part of $f(x) x^{n} + g(x)$ is irreducible
or identically one unless one of the following holds:
\begin{itemize}
\item[(i)]
The polynomial
$-f(x)g(x)$ is a $p^{\rm th}$ power for some  prime $p$ dividing $n$.
\item[(ii)]
For either $\varepsilon = 1$ or $\varepsilon = -1$, one
of $\varepsilon f(x)$ and $\varepsilon g(x)$ is a $4^{\rm th}$ power, the other
is $4$ times a $4^{\rm th}$ power, and $n$ is divisible by $4$.
\end{itemize}
\end{theorem}

Note that when (i) or (ii) hold in the above theorem, $f(x)x^n+g(x)$
is reducible by an apparent factorization.

\subsection{Proof of Theorem \ref{thm.whitehead}}
\label{sub.proof.thm.whitehead}

We set $f(x) = (x+1)^{m}$, $g(x) = -(x-1)^{m}$ and
\[
F_{k,m}(x) = f(x) x^{k} + g(x) = (x+1)^{m} x^{k} - (x-1)^{m},
\]
where $k$ is a positive integer.  
Observe that, since $m$ is odd, $\gcd(4n+m,m) = \gcd(n,m)$. 
Our goal is to show that if $k = 4n+m$ is sufficiently large and 
$\gcd(k,m) = 1$, then
$F_{k,m}(x)$ has $x^{2}+1$ as a factor and $F_{k,m}(x)/(x^{2}+1)$ is irreducible.  
We apply Theorem \ref{ydeg1} with $n$ replaced by $k$.  
For $k = 4n+m$ and $\gcd(k,m) = 1$, we show the following.
\begin{itemize}
\item
Each of (i) and (ii) in Theorem \ref{ydeg1} do not hold.
\item
The only cyclotomic factor of $F_{k,m}(x)$ is $x^{2}+1$, which has 
multiplicity $1$.
\item  
The polynomial $F_{k,m}(x)$ has no irreducible reciprocal non-cyclotomic factors.  
\end{itemize}
Theorem \ref{ydeg1} will then imply that $F_{k,m}(x)/(x^{2}+1)$ is 
irreducible for
$k = 4n+m$ sufficiently large and $\gcd(n,m) = 1$, completing the proof.

Suppose $k = 4n+m$ and $\gcd(k,m) = 1$.  
The polynomial $-f(x)g(x) = (x+1)^{m}(x-1)^{m}$ is a $p^{\rm th}$ power for
a prime $p$ if and only if $p|m$.  If $p|m$, then the condition $\gcd(k,m) = 1$ 
implies $p \nmid k$. Therefore, (i) with $n$ replaced by $k$, does not hold.
Furthermore, one sees that since $m$ is odd, the polynomials 
$\pm f(x) = \pm (x+1)^{m}$
and $\pm g(x) = \mp (x-1)^{m}$ cannot be a fourth power.  Thus, (ii) does 
not hold.  

Next, we show that the only cyclotomic factor of $F_{k,m}(x)$ is $x^{2}+1$.  
To do this, we first note
that if $w(x)$ is any irreducible reciprocal factor of $F_{k,m}(x)$, then 
$w(x)$ divides
\begin{align*}
x^{\deg F_{k,m}} F_{k,m}(1/x) 
&= x^{k+m} \bigg(  \bigg( \dfrac{1}x + 1 \bigg)^{m} \dfrac{1}{x^{k}} 
- \bigg(  \dfrac{1}x - 1 \bigg)^{m}  \bigg) \\[5pt]
&= - (1-x)^{m} x^{k} + (x+1)^{m}.
\end{align*}
Thus, $w(x)$ must be a divisor of
\[
(1-x)^{m} F_{k,m}(x) + (x+1)^{m} x^{\deg F_{k,m}} F_{k,m}(1/x) 
= - (1-x)^{m} (x-1)^{m} + (x+1)^{2m}.
\]
Since $m$ is odd, we deduce that $w(x)$ divides $(x+1)^{2m} + (x-1)^{2m}$.
  
Assume now that $\zeta = e^{2\pi i/s}$, where $s \in \mathbb Z^{+}$, 
is a root of $F_{k,m}(x)$.  
Since cyclotomic polynomials are reciprocal, we obtain that $\zeta$ 
must be a root of $(x+1)^{2m} + (x-1)^{2m}$.  Hence, $\zeta \ne 1$ and
\[
\bigg(  \dfrac{\zeta+1}{\zeta-1}  \bigg)^{2m} = -1.
\]
This implies that $(\zeta+1)/(\zeta-1)$ is a $(4m)^{\rm th}$ root of unity.  
Setting $\theta = 2\pi/s$, 
we see that
\[
\dfrac{\zeta+1}{\zeta-1} = \dfrac{e^{i\theta}+1}{e^{i\theta}-1}
= \dfrac{e^{i\theta/2}+e^{-i\theta/2}}{e^{i\theta/2}-e^{-i\theta/2}} 
= \dfrac{\cos (\theta/2)}{i \sin (\theta/2)}
= - i \cot (\theta/2).
\]
Since $(\zeta+1)/(\zeta-1)$ is a root of unity, this last expression is 
non-zero and necessarily 
purely imaginary.  
Also, $\theta = 2\pi/s \in (0,\pi]$. 
Therefore, $(\zeta+1)/(\zeta-1) = -i$, and we must have $\cot (\theta/2) = 1$.  
We see then that $s = 4$ and $\zeta = i$.  Observe that, since $k = 4n+m$,
we can deduce that
\[
\dfrac{F_{k,m}(\zeta)}{(\zeta-1)^{m}} 
=  \bigg(  \dfrac{\zeta+1}{\zeta-1}  \bigg)^{m} \zeta^{m} \zeta^{4n} - 1
= (-i)^{m} i^{m} i^{4n} - 1 = 1 - 1 = 0.
\]
Thus, $x^{2}+1$ is a factor of $F_{k,m}(x)$.  Since
\[
F'_{k,m}(x) = m (x+1)^{m-1} x^{k} + k (x+1)^{m} x^{k-1} - m (x-1)^{m-1},
\]
we also obtain
\begin{align*}
\dfrac{F'_{k,m}(\zeta)}{(\zeta-1)^{m-1}} &=  
m \bigg(  \dfrac{\zeta+1}{\zeta-1}  \bigg)^{m-1} \zeta^{m} \zeta^{4n} 
+ k \bigg(  \dfrac{\zeta+1}{\zeta-1}  \bigg)^{m-1} (\zeta+1) \zeta^{m-1} 
\zeta^{4n} - m \\[5pt]
&= m (-i)^{m-1} i^{m} i^{4n} + k (-i)^{m-1} (i+1) i^{m-1} i^{4n} - m 
= (k+m)i + (k-m) \ne 0.
\end{align*}
Thus, $(x^{2}+1)^{2}$ is not a factor of $F_{k,m}(x)$, and 
$F_{k,m}(x)/(x^{2}+1)$ has no cyclotomic factors.

We now assume that $F_{k,m}(x)$ has an irreducible 
reciprocal non-cyclotomic factor.
Calling this factor $w(x)$, we recall that $w(x)$ must divide 
$(x+1)^{2m} + (x-1)^{2m}$.  
Equivalently, we have
\begin{equation}\label{norecipeq1}
(x+1)^{2m} \equiv - (x-1)^{2m} \pmod{w(x)}.
\end{equation}
Since $w(x)$ divides $F_{k,m}(x)$, we also have
\begin{equation}\label{norecipeq2}
(x+1)^{m} x^{k} \equiv (x-1)^{m} \pmod{w(x)}.
\end{equation}
Squaring both sides of \eqref{norecipeq2} and substituting
from~\eqref{norecipeq1}, we see that
\[
- (x-1)^{2m} x^{2k} \equiv (x-1)^{2m} \pmod{w(x)}.
\]
Hence, the polynomial $w(x)$ divides $(x-1)^{2m} (x^{2k}+1)$.  As $w(x)$
is irreducible and non-cyclotomic and as all the irreducible factors of
$(x-1)^{2m} (x^{2k}+1)$ are cyclotomic, we obtain a contradiction.  Thus, the
polynomial $F_{k,m}(x)$ has no irreducible reciprocal non-cyclotomic factors.
The proof is complete.
\qed


\subsection{Proof of Theorem~\ref{thm.whitehead2}}

The main part of the theorem is a direct consequence of the work of Neumann-Reid
from~\cite[Thm.~6.2, Eqn.~(6.7)]{NR} which was mentioned in the introduction
combined with Theorem~\ref{thm.whitehead}. What is left to explain is the
stated recursion~\eqref{Tmnrec}.

We want to justify that 
\begin{equation}
\label{filgaradd1}
(x^{2}+1)  x^{2n+m-1} T_{m,n}\bigg(  x - \dfrac{1}{x}  \bigg) = F^{W}_{m,n}(x),
\end{equation}
where $T_{m,n}(z)$ is defined by the recursion~\eqref{Tmnrec}.
Recall $z = x - 1/x$. One checks directly that~\eqref{filgaradd1} holds for
$(m,n) \in \{ (0,0), (0,1), (1,0), (1,1) \}$.  Suppose~\eqref{filgaradd1} holds
when $(m,n)$ is replaced by $(m-1,n)$ and $(m-2,n)$. Then making use of the
recursion for $T_{m,n}(z)$ with $z=x-1/x$, we obtain

\begin{small}
\begin{align*}
 (x^{2}+1)  \,x^{2n+m-1} T_{m,n}(z) &  \\ 
 &\hspace{-3cm}=  (x^{2}+1)  \,x^{2n+m-1}  (z+2) \,
 T _{m-1,n}(z) 
 - (x^{2}+1)  \,x^{2n+m-1} \bigg(  x - \dfrac{1}{x}  \bigg) \,
 T _{m-2,n}(z) \\[8pt]
 &\hspace{-3cm}=  (x^{2}+1)  \,x^{2n+(m-1)-1}  \big(  x^{2} - 1 + 2x \big) \, T _{m-1,n}
 (z) - (x^{2}+1)  \,x^{2n+(m-2)-1} \big(  x^{3} - x  \big) \,T _{m-2,n}
 (z) \\[8pt]
 &\hspace{-3cm}=  \big(  x^{2} - 1 + 2x \big) \, F^{W}_{m-1,n}(x) - \big(  x^{3} - x  \big)
 \,F^{W}_{m-2,n}(x).
\end{align*}
\end{small}
Thus, we want to show
\[
  \big(  x^{2} - 1 + 2x \big) \, F^{W}_{m-1,n}(x) - \big(  x^{3} - x  \big)
  \,F^{W}_{m-2,n}(x) = F^{W}_{m,n}(x).
\]
This identity is easily verified from the definition of $F^{W}_{m,n}(x)$ on noting
that
\begin{gather*}
(x^2+2x-1)(x^2+x)-(x^3-x) = \big(x (x+1)\big)^{2} 
\shortintertext{and}
(x^2+2x-1)(x-1)-(x^3-x) = (x-1)^{2}.
\end{gather*}
The above shows that~\eqref{filgaradd1} holds if it holds for $(m,n)$ replaced
by $(m-1,n)$ and $(m-2,n)$. A similar argument shows that~\eqref{filgaradd1}
holds if it holds for $(m,n)$ replaced by $(m,n-1)$ and $(m,n-2)$. 
This is enough to deduce that~\eqref{filgaradd1} holds in general.
\qed


\section{Factoring fewnomials}
\label{sec.fewnomials}{}

\subsection{The cyclotomic part of $G^T_{a,b}(x)$}
\label{sub.cycloPab}
In this section, we prove Theorem~\ref{thm.teichmuller} (a).
Suppose $\zeta = \zeta_{n} = e^{2\pi i/n}$ is a root of $G^T_{a,b}(x)$, and hence 
\begin{equation}
  \label{teichmulleraeq1}
\zeta^{a+2b} +\zeta^{a} - \zeta^{2a+b} - \zeta^{b} - \zeta^{a+b} = 0.
\end{equation}
The strategy will be as follows.  
We will show that 
\[
  (*) \ \,\text{\parbox[t]{14cm}{there are positive integers
      $n_{0}$, $d_{0}$, $a_{0}$ and $b_{0}$ such that $n_{0} \le 12$, 
$n = n_{0} d_{0}$, $a = a_{0} d_{0}$ and $b = b_{0} d_{0}$.}}
\]
Then
\[
\zeta^{a} = \zeta_{n_{0}d_{0}}^{a_{0}d_{0}} = \zeta_{n_{0}}^{a_{0}} 
\quad \text{ and } \quad
\zeta^{b} = \zeta_{n_{0}d_{0}}^{b_{0}d_{0}} = \zeta_{n_{0}}^{b_{0}} 
\]
so that \eqref{teichmulleraeq1} implies
\begin{equation}\label{teichmulleraeq2}
  \zeta_{n_{0}}^{a_{0}+2b_{0}} +\zeta_{n_{0}}^{a_{0}} - \zeta_{n_{0}}^{2a_{0}+b_{0}}
  - \zeta_{n_{0}}^{b_{0}} - \zeta_{n_{0}}^{a_{0}+b_{0}} = 0.
\end{equation}
As $n_{0} \le 12$, to obtain all the solutions to \eqref{teichmulleraeq2}, 
we begin by calculating each possibility for $1 \le n_{0} \le 12$,
$1 \le a_{0} \le n_{0}$ and $1 \le b_{0} \le n_{0}$ for
which \eqref{teichmulleraeq2} holds. This can be accomplished by checking for
each such $n_{0}$, $a_{0}$ and $b_{0}$ whether
$x^{a_{0}+2b_{0}} +x^{a_{0}} - x^{2a_{0}+b_{0}} - x^{b_{0}} - x^{a_{0}+b_{0}}$
is divisible by $\Phi_{n_{0}}(x)$, as \eqref{teichmulleraeq2} will hold if
and only if this divisibility by $\Phi_{n_{0}}(x)$ holds. Once we determine the
set $\mathcal T$ of all triples $(n_{0}, a_{0}, b_{0})$ with
$1 \le n_{0} \le 12$, $1 \le a_{0} \le n_{0}$ and $1 \le b_{0} \le n_{0}$
satisfying \eqref{teichmulleraeq2}, we deduce that all solutions
to \eqref{teichmulleraeq2} with $n_{0} \le 12$ are given by triples
$(n_{0}, a_{0}, b_{0})$ of positive integers satisfying
\begin{equation}
  \label{teichmulleraeq3}
a_{0} \equiv a'_{0} {\hskip -4pt}\pmod{n_{0}}, \quad
b_{0} \equiv b'_{0} {\hskip -4pt}\pmod{n_{0}}, \quad \text{and} \quad
(n_{0},a'_{0},b'_{0}) \in \mathcal T.
\end{equation}
One can reverse these steps to see that \eqref{teichmulleraeq3}
implies \eqref{teichmulleraeq2} implies \eqref{teichmulleraeq1}, where
$n = n_{0} d_{0}$, $a = a_{0} d_{0}$ and $b = b_{0} d_{0}$ for some positive
integer $d_{0}$.  We deduce then that $G^T_{a,b}(x)$ is divisible by $\Phi_{n}(x)$
if and only if $n = n_{0} d_{0}$, $a = a_{0} d_{0}$ and $b = b_{0} d_{0}$
where $d_{0} \in \mathbb Z^{+}$ and \eqref{teichmulleraeq3} holds.
Following a direct computation of $\mathcal T$ as above, the divisibility
of $G^T_{a,b}(x)$ by $C^T_{a,b}(x)$ then follows.  What remains to be shown
for part (a) of Theorem~\ref{thm.teichmuller}, is that
if \eqref{teichmulleraeq1} holds, then ($*$) holds. Furthermore, we need to
justify that each of the irreducible polynomials $\Phi_{n}(x)$ that divides
$G^T_{a,b}(x)$ appears as a factor in $G^T_{a,b}(x)$ with multiplicity exactly one.

For $k$ a positive integer relatively prime to $n$, let $\phi_{k}$ be the
automorphism of $\mathbb Q(\zeta)$ fixing $\mathbb Q$ given by
$\phi_{k}(\zeta) = \zeta^{k}$. Let $d = \gcd(a,n)$, and define positive integers
$a'$ and $n'$ by $a = d a'$ and $n = d n'$. Then $\gcd(a',n') = 1$.  
Letting
\[
a'' = a' + n' \prod_{\substack{p \mid d \\ p \,\nmid \,a'}} p,
\]
one can see that $\gcd(a'',n) = 1$ by considering separately the three cases
where a prime divides $n'$, where a prime divides $d$ but neither $n'$ nor $a'$,
and where a prime divides both $d$ and $a'$. Note that in the latter case,
since $\gcd(a',n') = 1$, we have the prime does not divide $n'$. Since
$\gcd(a'',n) = 1$, there exists a $k \in \mathbb Z$ such that
$k a'' \equiv 1 \pmod{n}$. Note that $k a'' \equiv 1 \pmod{n}$ implies
$\gcd(k,n) = 1$.  
Since $a''$ and $a'$ differ by a multiple of $n'$, we obtain
\[
\zeta_{n}^{ka} = \zeta_{n'}^{ka'} = \zeta_{n'}^{ka''} = \zeta_{n'}.
\]
We deduce from $\phi_{k}\big(P^{T}(\zeta^{a},\zeta^{b})\big) = 0$ that
\begin{equation}
  \label{thm3pt1(a)eq1}
  \zeta^{kb} + \zeta^{-kb} = \zeta_{n}^{ka} + \zeta_{n}^{-ka} + 1
  = \zeta_{n'} + \zeta_{n'}^{-1} + 1.
\end{equation}
In particular, if $n' > 6$, then 
\[
2 \ge |\zeta^{kb} + \zeta^{-kb}| = |\zeta_{n'} + \zeta_{n'}^{-1} + 1| 
= \cos(2\pi/n') + \cos(-2\pi/n') + 1 > 2,
\]
which is impossible.  Thus, $n' \le 6$.  Recall $\gcd(k,n) = 1$.  Let
$k' \in \mathbb Z$ satisfy $k k' \equiv 1 \pmod{n}$. 

If $n' = 6$, then \eqref{thm3pt1(a)eq1} implies $\zeta^{kb} + \zeta^{-kb} = 2$.
Thus, $\zeta^{kb} = 1$. Applying $\phi_{k'}$, we deduce that $\zeta^{b} = 1$.
As $\zeta_{n}^{b} = \zeta_{n'd}^{b} = \zeta_{6d}^{b}$, we deduce that $b$ is a
multiple of $6d$. Thus, at this point we have $n = n'd = 6d$, $a = a'd$ and
$b = b'd$ for some positive integer $b'$ divisible by $6$.  
Taking $d_{0} = d$ and $n_{0} = 6$, we see that ($*$) holds in this case.

If $n' = 5$, then  \eqref{thm3pt1(a)eq1} implies 
\[
\zeta^{kb} + \zeta^{-kb} = \zeta_{5} + \zeta_{5}^{-1} + 1 = \dfrac{1+\sqrt{5}}{2}.
\]
Letting $\alpha_{1} = (1+\sqrt{5})/2$ and $\alpha_{2} = (1-\sqrt{5})/2$, we see
that $\zeta^{kb}$ is a root of
\[
  (x^{2} - \alpha_{1} x + 1) (x^{2} - \alpha_{2} x + 1)
  = x^{4} - x^{3} + x^{2} -  x + 1 = \Phi_{10}(x).
\]
Thus, $\Phi_{10}(\zeta^{kb}) = 0$, and applying $\phi_{k'}$, we see that
$\Phi_{10}(\zeta^{b}) = 0$. We deduce that $n/\gcd(n,b) = 10$ so that
$n = 10 d'$ and $b = b'd'$ for some positive integers $d'$ and $b'$.  
Also, $n = n'd = 5d$ implies $d = 2d'$ so that $a = a'd = 2a' d'$.  
Taking $d_{0} = d'$ and $n_{0} = 10$, we see that ($*$) holds in this case.

If $n' = 4$, then using \eqref{thm3pt1(a)eq1} and applying $\phi_{k'}$ gives
$\zeta^{b} + \zeta^{-b} = 1$, so $\zeta^{b}$ is a root of
$x^{2} - x + 1 = \Phi_{6}(x)$.  Thus, $n/\gcd(b,n) = 6$, so there are positive
integers $d'$ and $b'$ such that $n = 6d'$ and $b=b'd'$. Since $n = n'd = 4d$,
we also see that there is a positive integer $d''$ such that 
\[
d' = 2d'', \quad
n = 12d'', \quad
d = 3d'', \quad
a = a'd = 3a'd'', \quad \text{and} \quad
b = b'd' = 2b'd''.
\]
Therefore, with $d_{0} = d''$ and $n_{0} = 12$, we obtain ($*$) in this case.

If $n' = 3$, then $3$ divides $n$. 
In this case, using \eqref{thm3pt1(a)eq1} and applying $\phi_{k'}$ leads to
$\zeta^{b} + \zeta^{-b} = 0$, so $\zeta^{b} = \pm i$. We deduce $n/\gcd(b,n) = 4$.
Since $3$ divides $n$, we can write $n = 12 d'$ for some positive integer $d'$.  
Then $\gcd(b,n) = 3d'$ so that $b = 3b'd'$ for some positive integer $b'$.  Since
$n = n'd = 4d$, we have $d = 3d'$. Hence, $a = a'd = 3a'd'$.  Thus, using
$d_{0} = d'$ and $n_{0} = 12$, we see ($*$) holds in this case.

If $n' = 2$, then $2$ divides $n$.
Here, using \eqref{thm3pt1(a)eq1} and applying $\phi_{k'}$ gives
$\zeta^{b} + \zeta^{-b} = -1$, so $\zeta_{b}$ is a root of $x^{2}+x+1 = \Phi_{3}(x)$.  
Thus, $n/\gcd(b,n) = 3$, so $3$ also divides $n$.  We write $n = 6d'$.  
Then 
\[
\gcd(b,n) = 2d', \quad
b = 2b'd', \quad
d = 3d', \quad \text{and} \quad
a = 3a'd',
\]
where $a'$ and $b'$ are positive integers.  
With $n_{0} = 6$ and $d_{0} = d'$, we see that  ($*$) holds in this case as well.

If $n' = 1$, then \eqref{thm3pt1(a)eq1} implies $\zeta^{kb} + \zeta^{-kb} = 3$,
which is impossible since the left-hand side is at most $2$.  So the case
$n' = 1$ cannot occur. This finishes the proof that \eqref{teichmulleraeq1}
implies ($*$). 

To finish the proof of Theorem~\ref{thm.teichmuller} (a), we still need to
address the multiplicity of the cyclotomic factors.  
Assume now that $\zeta = \zeta_{n}$ is a root with multiplicity $> 1$ for
$G^{T}_{a,b}(x)$. Setting $H(x) = H_{a,b}(x)$ to be the numerator in \eqref{GTdef},
we deduce $H(\zeta) = H'(\zeta) = 0$. Since also $\zeta H'(\zeta) = 0$, we deduce
that $\zeta$ is a root of
\begin{equation}
  \label{thm3pt1(a)eq3}
U(x) := (a+2b) x^{a+2b} + a x^{a} -  (2a+b)x^{2a+b} - bx^{b} - (a+b)x^{a+b}.
\end{equation}
By what we have shown above,
we know that either (I), (II) or (III) holds.  

In the case of (I), we have $n = 6d_{1}$ so that $\zeta^{d_{1}} = \zeta_{6}$.  
We deduce $\zeta^{a} = \zeta^{a_{1}d_{1}} = \zeta_{6}^{a_{1}}$ and, similarly,
$\zeta^{b} = \zeta_{6}^{b_{1}}$. Now, it is simply a matter of substituting in
the four possible values of $(a_{1},b_{1})$ modulo $6$ for $(a,b)$ into
the \textit{exponents} of $U(x)$ in \eqref{thm3pt1(a)eq3} and verify that $U(x)$
does not have $\Phi_{6}(x)$ as a factor, giving a contradiction.
For example, if $(a_{1}, b_{1}) \equiv (3,2) \pmod{6}$, then 
\[
U(x) = (a+2b) x^{7} + a x^{3} -  (2a+b) x^{8} - b x^{2} - (a+b) x^{5} = 0.
\]
The remainder when $U(x)$ is divided by $\Phi_{6}(x) = x^{2} - x + 1$ is
\[
(a+2b) x - a -  (2a+b) (x-1) - b (x-1)  + (a+b) (x-1) = b (x+1).
\]
As the remainder must be identically $0$ for $\Phi_{6}(x)$ to divide $U(x)$,
we deduce $b = 0$, contradicting that $b \in \mathbb Z^{+}$. A very similar
argument works for each of the four possibilities for $(a_{1},b_{1})$ modulo
$6$ in (I); more precisely, in each case, the coefficient of $x$ in the remainder
being $0$ implies that either $a$ or $b$ is $0$, giving a contradiction.

In the case of (II), we obtain $\zeta^{d_{2}} = \zeta_{10}$ so that
$\zeta^{a} = \zeta_{10}^{a_{2}}$ and $\zeta^{b} = \zeta_{10}^{b_{2}}$. One similarly
checks each of the remainders obtained by substituting the eight possibilities
for $(a_{2},b_{2})$ for $(a,b)$ in the exponents of $U(x)$ and dividing by
$\Phi_{10}(x)$. As $\deg \Phi_{10}(x) = 4$, the remainders will have degree at
most $3$ with coefficients possibly depending on $a$ and $b$. In each case,
the coefficient of $x^{3}$ or the coefficient of $x^{2}$ being $0$ implies
that $a$ or $b$ is $0$, giving a contradiction.  For example, in the case
$(a_{2},b_{2}) = (4, -3) \pmod{10}$, the remainder of $U(x)$ divided by
$\Phi_{10}(x)$ is $-2b x^3 + (b-a)x^2 - bx + a + b$, so a contradiction is
obtained since the coefficient of $x^{3}$ being $0$ implies $b = 0$.

In the case of (III), we deduce  $\zeta^{d_{3}} = \zeta_{12}$ so that
$\zeta^{a} = \zeta_{12}^{a_{3}}$ and $\zeta^{b} = \zeta_{12}^{b_{3}}$.  
We consider the remainders obtained by substituting each of the eight
possibilities for $(a_{3},b_{3})$ modulo $12$ for $(a,b)$ in the exponents of
$U(x)$ and dividing by $\Phi_{12}(x)$. Note that $\deg \Phi_{12}(x) = 4$.
In each case, the coefficient of $x^{3}$ in the remainder is non-zero if $a$
and $b$ are non-zero, so we obtain a contradiction.

Theorem~\ref{thm.teichmuller} (a) follows.
\qed

\vskip 5pt \noindent
\textbf{Comment:}  We opted for a self-contained argument above, in part because
of the simplicity of the argument. There is a nice general algorithm for
establishing the cyclotomic factors of polynomials given by Granville and
Rudnick~\cite[Section 3]{Granville}. This leads to narrowing down consideration
to $m^{\rm th}$ roots of unity where $m$ divides $30$. The multiplicity of the roots
would also need to be considered separately as done here. Another approach would
be to take advantage of the fact that replacing $x$ and $y$ with cyclotomic numbers
in \eqref{eq.teichmuller} allows one to rewrite $P^{T}(x,y) = 0$ as a sum of
three cosine values equal to $0$, where one depends on $x$, one depends on $y$ and
one is $\cos(2\pi/3)$.  Then a computational step in the recent interesting work
by Kedlaya, Kolpakov, Poonen and Rubinstein \cite[Theorem~5.1]{KKPR} can
be applied, where again the multiplicity of the roots would still need to be
considered.

\subsection{The non-reciprocal part of $G^T(x^a,x^b)$}

In this section, we establish Theorem~\ref{thm.teichmuller}~(b). 
For distinct $a$ and $b$, one checks that the exponents appearing in the numerator
of $G^T(x^a,x^b)$ in the last expression in \eqref{GTdef} are distinct;
in other words, $G^T(x^a,x^b)$ consists of exactly five terms with coefficients
$\pm 1$. Furthermore, the polynomial $G^T(x^a,x^b)$ is reciprocal.  A reciprocal
polynomial cannot have exactly one irreducible non-reciprocal factor, so the
non-reciprocal part of $G^T(x^a,x^b)$ is either reducible or identically $1$.  
To establish Theorem~\ref{thm.teichmuller}~(b), we justify that the non-reciprocal
part of $G^T(x^a,x^b)$ is not reducible.

A recent paper by Murphy, Vincent and the first author \cite[Theorem 1.2]{FMV}
provides a description of every polynomial $f(x)$ consisting of five terms and
coefficients from $\{ \pm 1 \}$ that has a reducible non-reciprocal part.  Although
this reference is ideal for our situation, the result in \cite{FMV} needs to be
unraveled as it indicates that such an $f(x)$ must come from a ``modification" of
one of $13$ possible explicit classifications.  The details of how to interpret
the result are in \cite{FMV}.  We give some explanation here as well by using as
a concrete example one of the more complicated of the $13$ classifications for
such an $f(x)$. 

Let 
\begin{equation}
  \label{FMVeg}
F(x,y) = x^5y^2-x^3y^2-x^3y-xy-1 = ( x^{3}y+x^{2}y+1 )  ( x^{2}y-xy-1 ).
\end{equation}
Then one of the 13 classifications noted above from \cite{FMV} is that $f(x)$ has
a reducible non-reciprocal part if $f(x)$ is of the form $G(x^{t},x^{u})$ where
$G(x,y)$ is a modification of the polynomial $F(x,y)$ and $t$ and $u$ are positive
integers. We first describe what these modifications look like by indicating how
they are constructed. For a general polynomial $H(x,y) \in \mathbb Z[x,y]$, define
\[
\widetilde{H}(x,y) = x^{k} y^{\ell} H(1/x,1/y),
\]
where $k$ and $\ell$ are integers chosen as small as possible so that
$\widetilde{H}(x,y) \in \mathbb Z[x,y]$. If $H(x,y) = \pm \widetilde{H}(x,y)$,
we say that $H(x,y)$ is reciprocal. Set $S_{1} = \{ \pm F(x,y) \}$.  Given $S_{j}$
with $j \ge 1$, we construct $S_{j+1}$ as follows. For each $F_{0}(x,y) \in S_{j}$,
we begin by factoring it in $\mathbb Z[x,y]$ as $F_{0}(x,y) = F_{1}(x,y) F_{2}(x,y)$,
where $ F_{1}(x,y)$ and $F_{2}(x,y) $ are irreducible non-reciprocal polynomials. 
A check can be done, as one proceeds, to see that starting with any $F(x,y)$
from Theorem~1.2 in \cite{FMV}, each $F_{0}(x,y) \in S_{j}$ for every $j$ can be
factored in this way. In other words, for the purposes of Theorem~1.2 in \cite{FMV},
which we are using, each $F_{0}(x,y)$ in each $S_{j}$ always factors as a product of 
exactly two irreducible non-reciprocal polynomials.   
The ordering of these two factors is irrelevant but to be fixed for each
$F_{0}(x,y) \in S_{j}$ in what follows.
We set
\[
F_{0}^{*}(x,y) = F_{1}(x,y) \widetilde{F}_{2}(x,y)
\]
and
\[
S_{j+1} = S_{j} 
\bigcup \bigg(  \bigcup_{F_{0} \in S_{j}}
\big\{ \pm \widetilde{F}_{0}(x,y) \big\}  \bigg)
\bigcup \bigg(  \bigcup_{F_{0} \in S_{j}}
\big\{ \pm F_{0}^{*}(x,y), \pm F_{0}(-x,y), \pm F_{0}(x,-y) \big\}  \bigg).
\]
A key point is that there will be a $J \in \mathbb Z^{+}$ such that $S_{J} = S_{J+1}$
and as a consequence we obtain $S_{j} = S_{J}$ for all $j \ge J$.  
We construct the sets $S_{j}$ until such a $J$ occurs.  Then the set $S_{J}$
corresponds to the set of all modifications of $F(x,y)$.  
For $F(x,y)$ defined by \eqref{FMVeg}, we obtain $J = 5$ and $|S_{J}| = 32$.

Next, we need to determine whether there are choices for $a$ and $b$ for which
$G^T(x^a,x^b)$ is of the form given by $G(x^{t},x^{u})$ for each $G(x,y) \in S_{J}$.
We simplified the process by only looking at such $G(x,y)$ with two positive
coefficients and three negative coefficients, as in $G^T(x^a,x^b)$.  
For $F(x,y)$ defined by \eqref{FMVeg}, there were $8$ such polynomials in $S_{J}$.  
Fix $G(x,y)$ as one of these. We then want to compare the degrees of the positive
terms in $G(x^{t},x^{u})$ with the positive terms in $G^T(x^a,x^b)$ and the negative
terms in  $G(x^{t},x^{u})$ with the negative terms in $G^T(x^a,x^b)$. In the way of
an explicit example, for $F(x,y)$ as in \eqref{FMVeg}, one of the elements of
$S_{J}$ is $G(x,y) = x^5y^2 - x^4y - x^2y - x^2 + 1$.  Observe that
$G(x^{t},x^{u}) = x^{5t+2u} - x^{4t+u} - x^{2t+u}  - x^{2t} + 1$. We set
$c = -\min\{ a,b \}$ to be an unknown in the definition of $G^T(x^a,x^b)$ in the
last expression in \eqref{GTdef}. In order for $G^T(x^a,x^b)$ to be of the form
$G(x^{t},x^{u})$, then it is necessary and sufficient for 
\begin{gather*}
\{ a+2b+c, a+c \}
= \{ 5t+2u, 0 \} 
\intertext{and} 
\{ 2a+b+c, b+c, a+b+c \}
= \{ 4t+u, 2t+u, 2t \}.
\end{gather*}
Although it is easy in this case to compare sizes of the elements to determine,
for example, that $a+2b+c = 5t+2u$ and $a+c = 0$, we simply made all possible
equations from the equality of these sets, allowing $2$ possibilities for the
first set equality and $6$ possibilities for the second set equality. Thus, for
this choice of $F(x,y)$ and $G(x,y)$, we had 12 systems of equations to solve in
integers $a$, $b$, $c$, $t$ and $u$. We used Maple 2019.2 to run through all of
the choices for $F(x,y)$ and $G(x,y)$ as described above, and none of resulting
systems of equations had solutions with non-zero $a$ and $b$.  
Thus, Theorem~\ref{thm.teichmuller}~(b) follows.

\subsection{Proof of Theorem~\ref{thm.teichmuller}}

We have already proven parts (a) and (b) of the theorem. For the last
part, we need to assume Lehmer's Conjecture~\cite{lehmer}.
Recall the Mahler measure of a polynomial $p(x_1,\dots, x_r)$ in $r$ variables
is given by
\begin{equation*}
M(p(x_1,\dots, x_r))=
\exp \bigg( \int_{0}^{1} \int_{0}^{1} \cdots \int_{0}^{1} 
\log |p(e^{2\pi i t_{1}},e^{2\pi i t_{2}},\ldots,e^{2\pi i t_{r}})| 
\,d t_{1} \,d t_{2} \ldots \,d t_{r} \bigg).
\end{equation*}
For a vector $\vec{s} = \langle s_{1},s_{2},\ldots,s_{r} \rangle$ of integers,
we define $H(\vec{s}\,) = \max_{1 \le j \le r} \{ |s_{j}| \}$. Now, given a vector
$\vec{a} = \langle a_{1},a_{2},\ldots,a_{r} \rangle$ of positive integers,
we consider all non-zero vectors $\vec{s} = \langle s_{1},s_{2},\ldots,s_{r} \rangle$
of integers perpendicular to $\vec{a}$, so satisfying $\sum_{j=1}^{r} s_{j} a_{j} = 0$, 
and define $q(\vec{a}\,)$ to be the minimum value of $H(\vec{s}\,)$ over all
such $\vec{s}$. Lawton~\cite[Theorem~2]{Lawton} has shown
\begin{equation*}
\lim_{q(\vec{a}\,) \to \infty}
M(p(x^{a_1},\dots, x^{a_r})) = M(p(x_1,\dots, x_r)).
\end{equation*}

Lehmer's Conjecture asserts that
\begin{gather*}
\label{lehmer}
M(p(x)) \geq c=1.176280818\ldots, 
\shortintertext{where}
c=M(x^{10}+x^9-x^7-x^6-x^5-x^4-x^3+x+1).
\end{gather*}
A computation gives $M(xyP^T(x,y))=1.28573\ldots$. 
We only require this numerical estimate, but we note that, as a consequence of
work of Rogers and Zudilin \cite{RogersZudilin}, an explicit expression of the
value of $M(xyP^T(x,y))$ can be given in terms of an $L$-series of a conductor
$15$ elliptic curve. Since $C^T_{a,b}(x)$ is the cyclotomic factor of
$G^T_{a,b}(x)$, Lawton's theorem implies that for $a/b$ or $b/a$ sufficiently
large (so that $q(\langle a, b \rangle)$ is large), one has
$M(G^T_{a,b}(x)/C^T_{a,b}(x)) = M(x^{a+b}P^T(x^{a},x^{b}))$ is near $1.28573\dots$
and hence $< c^{2}$. Mahler measure behaves multiplicatively on a product
of polynomials. Thus, if $G^T_{a,b}(x)/C^T_{a,b}(x)$ were reducible, it
would have a factor with Mahler measure less than $c$ violating
Lehmer's Conjecture. This completes the proof of Theorem~\ref{thm.teichmuller}~(c).
\qed


\appendix

\section{Dehn fillings of the Whitehead link}
\label{sec.dehn}

In this section, which is not needed for the results of the paper,
we recall the computation of the trace field of Dehn fillings 
of the Whitehead link, due to Neumann-Reid~\cite{NR}. The computation
works for every cusped hyperbolic manifold and uses standard material on ideal
triangulations, gluing equations and character varieties.

\subsection{Ideal triangulations, gluing equations and character 
varieties}
\label{sub.ideal}

An ideal tetrahedron is the convex hull of four points in the boundary
$\pt \BH^3$ of 3-dimensional hyperbolic space $\BH^3$~\cite{Th}. In the upper 
half-space model $\BH=\BC \times (0,\infty)$, the boundary can be identified
with the 2-dimensional sphere $\BC \cup \{\infty\}$. Setting three of the
four points of an ideal tetrahedron to $0$, $1$ and $\infty$ allows the
last point to be a complex number $z \in \BC\setminus\{0,1\}$. This number
is the {\em shape} of the ideal tetrahedron, equal to a cross ratio of
four points in $\BC$, and it is well-defined up to a six-fold orbit ambiguity.

Ideal tetrahedra were used by Thurston as building blocks of cusped hyperbolic 
3-manifolds as follows. Start from a topological ideal triangulation $\calT$ 
of an
oriented 3-manifold $M$ with torus boundary components (such a triangulation
always exists) and choose shape parameters $z_1,\dots,z_N$, one for each
ideal tetrahedron. Then, one can write down a system of gluing equations
around each edge of the ideal triangulation. Each equation 
implies that a Laurent monomial in $z_i$ and $w_i=1-z_i$ is equal to $\pm 1$. The
solution to the system of equations is an affine variety $G(\calT)$
in $(\BC\setminus\{0,1\})^N$. If $M$ has $c$ cusps, then typically 
$G(\calT)$ is a $c$-dimensional variety, although it is possible that 
$G(\calT)$
is empty. At any rate, each solution to the 
gluing equations gives via the developing map 
a representation of $\pi_1(M)$ (the fundamental group 
of $M$) to $\PSL(2,\BC)$, the isometry group of $\BH^3$.

If in addition one imposes a completeness equation around each cusp of $M$
(these equations also have the same shape as the gluing equations of the
edges of $\calT$), and finds a solution of the logarithmic system of equations,
this obtains a concrete description of the hyperbolic structure of $M$ in
terms of the shapes of $\calT$. In this case, the set of edge and
cusp equations is zero-dimensional and the shapes are algebraic numbers.
 
This method of describing hyperbolic manifolds may sound abstract. However,
Weeks developed a program {\tt SnapPea} that uses ideal triangulations and
numerically finds a hyperbolic structure. A modern version of 
the program, {\tt SnapPy}, is available from~\cite{snappy}. Among other
things, {\tt SnapPy} gives certified exact arithmetic computations.
For more information on the combinatorics of ideal triangulations of 
3-manifolds, see~\cite{NZ}. 

\subsection{Fillings of the Whitehead link}
\label{sub.whitehead}
 
In this section, we review the results from~\cite{NR} 
regarding fillings of the complement of the Whitehead link 
(shown on the left of 
Figure~\ref{f.2links}). The starting point is Thurston's observation that the
complement of the Whitehead link $W$ is a hyperbolic manifold obtained by 
face-pairings of a regular ideal octahedron~\cite[Fig.10]{NR}. The ideal
octahedron can be decomposed into four ideal tetrahedra, each of shape 
$\sqrt{-1}$~\cite[Fig.11]{NR}. This gives an ideal triangulation $\calT$
of $W$ with four tetrahedra. 
The gluing equations variety $G(\calT)$ is an affine surface.
Neumann-Reid~\cite[Sec.6.2]{NR} prove that the subset $G(\calT)^{\pr}$ of 
$G(\calT)$ which consists of solutions that are complete at the red cusp 
(of Figure~\ref{f.2links}) is a rational curve parametrized by
$$
\BC^* \longto G(\calT)^{\pr}, \qquad x \mapsto 
(z_1,z_2,z_3,z_4)=(x,-x^{-1},x,-x^{-1})
$$
The curve $G(\calT)^{\pr}$ contains the special point $x=i$ corresponding to
the hyperbolic structure on $W$. The (square of the) holonomy around the 
meridian and longitude of the blue cusp are rational functions $\mu(x)$ and
$\lambda(x)$ of $x$ on the curve $G(\calT)^{\pr}$. If $W(m/n)$ denotes 
the Dehn filling on the blue cusp of the Whitehead link $W$, then one adds the
equation
\begin{equation}
\label{eq.ml}
\mu(x)^m \lambda(x)^n =1
\end{equation}
Since $G(\calT)^{\pr}$ is a curve, if $\mu$ is a nonconstant function,
it follows that $\mu(x)^{-m}$ and $t=\lambda(x)$ are polynomially dependent, so
$p_m(t,\mu(x)^{-m})=0$ for some polynomial $p_{m}$. This brings \eqref{eq.ml}
to the form 
\begin{equation}
\label{eq.xl}
p_m(t,t^n) = 0
\end{equation}
as stated in the introduction. In the particular case of the Whitehead link,
Neumann-Reid~\cite[Eqn.~6.7]{NR} prove that \eqref{eq.ml} is simply
\[
(x(x+1))^m x^{4n}-(x-1)^m =0.
\]
Finally, to describe the trace field of $W(m/n)$, Neumann-Reid~\cite[p.~299]{NR} 
use the face-pairing of the ideal octahedron (4 pairs of faces, identified 
pairwise) and compute explicitly the four corresponding matrices of 
$\PSL(2,\BC)$. The entries of these matrices are polynomials in $x-x^{-1}$,
and the trace field that they generate is $\BQ(x-x^{-1})$.

We furthermore should point out a theorem of Hodgson-Meyerhoff-Weeks 
concerning pairs of {\em scissors congruent} hyperbolic manifolds 
$(M,M')$, that is pairs of non-isometric hyperbolic manifolds that are 
obtained by the same set of ideal tetrahedra, assembled together in a 
combinatorially different manner.

\begin{theorem}
\label{thm.HMW}\cite{HMW}
For coprime integers $(4m,n)$, the pair 
\[
\big(W(4m/n),-W(4m/(-n-2m)\big)
\] 
is geometrically similar and has a common 2-fold cover.
\end{theorem}
Here, $-M$ denotes the orientation reversed manifold of $M$.
 

\section*{Acknowledgments}

S.G. wishes thank BoGwang Jeon for enlightening conversations.


\bibliographystyle{hamsalpha}
\bibliography{biblio}
\end{document}